%%%%%%%%%%%%%%%%%%%
\input amstex
\documentstyle{amsppt}
\nologo
\NoBlackBoxes
%\NoRunningHeads
\mag=1200
%\pagewidth{5.2in}
%\pageheight{7.5in}
\hsize=31 pc  
\vsize=44 pc  
\hcorrection{5mm}

\topmatter
\title Knots and links without parallel tangents 
\endtitle
\author  Ying-Qing Wu$^1$
\endauthor
\leftheadtext{Y.-Q.~Wu}
\rightheadtext{Knots and links without parallel tangents}
\address Department of Mathematics, University of Iowa, Iowa City, IA
52242
\endaddress
\email  wu\@math.uiowa.edu
\endemail
%\keywords 
%\endkeywords
%\subjclass  Primary 57M25
%\endsubjclass
\thanks  $^1$ Partially supported by NSF grant \#DMS 9802558
\endthanks
%\abstract{}
%\endabstract
\endtopmatter

\document

\define\proof{\demo{Proof}}
\define\endproof{\qed \enddemo}
\define\a{\alpha}
\redefine\b{\beta}
\redefine\d{\delta}
\define\r{\gamma}
\redefine\R{{\Bbb R}}
\redefine\e{\epsilon}

\redefine\bdd{\partial}
\define\Int{\text{\rm Int}}

\baselineskip 15pt
\input epsf.tex

Steinhaus conjectured that every closed oriented $C^1$-curve has
a pair of anti-parallel tangents.  The conjecture is not true.  Porter
[Po] showed that there exists an unknotted curve which has no
anti-parallel tangents.  Colin Adams rised the question of whether
there exists a nontrivial knot in $\R^3$ which has no parallel or
antiparallel tangents.  In this paper we will solve this problem,
showing that any (smooth or polygonal) link $L$ in $\R^3$ is isotopic
to a smooth link $\hat L$ which has no parallel or antiparallel
tangents.  If ${\Cal S}(L)$ is the set of all smooth links isotopic to
$L$, then the subset $\hat {\Cal L}(L)$ of all $\hat L$ which has no
parallel or antiparallel tangents is not dense in ${\Cal S}(L)$ if it
is endowed with $C^{\infty}$ topology.  However, $\hat {\Cal L}(L)$ is
dense in $\Cal S(L)$ under $C^0$ topology.  We will show that any
neighborhood of $L$ contains such a link $\hat L$.  See Theorem 7
below.  The result has some impact on studying supercrossing numbers,
see the recent work of Pahk [Pa].

We refer the readers to [Ro] for concepts about knots and links.
Throughout this paper, we will use $I$ to denote a closed interval on
${\Bbb R}$.  Denote by $S^2$ the unit sphere in $\R^3$, and by $S_1$
the circle $S^2 \cap \R_{xy}$ on $S^2$, where $\R_{xy}$ denotes the
$xy$-plane in $\R^3$.  Denote by $Z[z_1, z_2]$ the set $\{ v=(x,y,z)
\in \R^3 \, | \, z_1\leq z \leq z_2\}$.  Similarly for $Y[y_1,
\infty)$ etc.  A curve $\b: I \to \R^3$ is an {\it unknotted curve\/}
in $Z[z_1,z_2]$ if (i) $\b$ is a properly embedded arc in
$Z[z_1,z_2]$, with endpoints on different components of $\bdd
Z[z_1,z_2]$, and (ii) $\b$ is rel $\bdd$ isotopic in $Z[z_1,z_2]$ to a
straight arc.

Given a curve $\a: I=[a,b] \to S^2$ and a positive function $f: I \to
\R_+ = \{x \in \R | x > 0\}$, we use $\b = \b(f, \a, t_0, v_0)$ to
denote the integral curve of $f\a$ with $\b(t_0) = v_0$, where $t \in
I$.  More explicitly, define
$$\b (t) = \b(f,\a,t_0,v_0)(t) = v_0 + \int_{t_0}^t f(t) \a(t) \,
dt.$$
When $t_0 = a$ and $v_0 = 0$, simiply denote it by $\b(f,\a)$.  

If $\r: [a,b] \to \R^3$ is a map and $[c,d]$ is a subinterval of
$[a,b]$, denote by $\r[a,b]$ the restriction of $\r$ on $[c,d]$.  If
$u, v$ are points in $\R^3$, denote by $e(u,v)$ the line segment with
endpoints at $u$ and $v$, oriented from $u$ to $v$.  Denote by
$d(u,v)$ the distance between $u$ and $v$.  Denote by $||e||$ the
length of $e$ if $e$ is a line segment or a vector in $\R^3$.  Thus
$d(u,v) = || e(u,v)|| = || u-v||$.

Given $n$ points $v_1, ..., v_n$ in $\R^3$, let $e_i = e(v_i,
v_{i+1})$.  The subscripts are always mod $n$ numbers.  Thus $e_n =
e(v_n, v_1)$.  In generic case, the union of these edges forms a knot,
denoted by $K(v_1, ..., v_n)$.  To avoid trivial case, we will always
assume $n \geq 4$.

\proclaim{Lemma 1} Let $K= K(v_1, ..., v_n)$ be a polygonal knot, and
let $N$ be a regular neighborhood of $K$.  Then there is a number $r>0$ such
that 

(i) $N$ contains the $r$-neighborhood $N(K)$ of $K$;

(ii) $K' = K(v'_1, v_1, ..., v'_n, v_n)$ is isotopic to $K$ in $N$ if
$d(v_i, v'_i) < r$; and 

(iii) $K'' = K(v'_1, ..., v'_n)$ is isotopic to $K$ in $N$ if $d(v_i,
v'_i) < r$.  \endproclaim

\proof Choose $r>0$ to satisfy (i) and $r < d/4$, where $d$ is the
minimal distance between non-adjacent edges of $K$.  The knot $K'$ is
contained in $N(K)$.  Let $D_i$ be the meridian disk of the
$r$-neighborhood $N(e_i)$ of the $i$-th edge $e_i = e(v_i, v_{i+1})$
of $K$, intersecting $e_i$ perpendicularly at its middle point $m_i$.
It is easy to check that the distance from $m_i$ to any edge $e_j$
($j\neq i$) is at least $d/2 = 2r$, hence $D_i$ is a meridian disk of
$N(K)$.  The edge $e(v'_i, v_i)$ is contained in an $r$-neighborhood
of $v_i$, hence is disjoint from all $D_j$.  Thus $D_i$ intersects
$K'$ at a single point on the edge $e(v_i, v'_{i+1})$, so the disks
$D_1, ..., D_n$ cut $N(K)$ into balls $B_1, ..., B_n$, each
intersecting $K'$ in an arc consisting of three edges, hence
unknotted.  Therefore $K'$ is isotopic to $K$ in $N(K)$.  This proves
(ii).

By (ii), both $K$ and $K(v'_1, ..., v'_n)$ are isotopic to $K(v_1,
v'_1, ..., v_n, v'_n)$ in $N$.  Therefore, they are isotopic to each
other, and (iii) follows.
\endproof

Denote by $C(\theta, u, v)$ the solid cone based at $u$ (the vertex of
the cone), open in the direction of $v$, with angle $\theta$.  More
explicitly, if we set up the coordinate system with $u$ the origin and
$v$ in the direction of $(0,0,1)$, then
$$C(\theta, u,v) = \{(x,y,z) \in \R^3 \, | \, z \geq \cot \theta
\sqrt{ x^2 + y^2} \}.$$

A smooth curve $\b: [a,b] \to B$ in a ball $B$ is {\it
$\theta$-allowable\/} if (i) $\b$ is properly embedded and unknotted
in $B$, (ii) the cones $C_a = C(\theta, \b(a), -\b'(a))$ and $C_b =
C(\theta, \b(b), \b'(b))$ are mutually disjoint, each intersecting $B$
only at its cone point.  

A smooth arc $\b : [a,b] \to \R^3$ is called an {\it
$\e$-suspension\/} if it is an embedding into an equilateral triangle
$\Delta$ in $\R^3$ with base the line segment $e = e(\b(a), \b(b))$
and height $\e$.  It is called a {\it round\/} $\e$-suspension if
furthermore it is a subarc of a round circle in $\R^3$, and
$||\b'(t)||$ is a constant function.  The line segment $e$ is called
the {\it base arc\/} of $\b$, and the disk bounded by $\b$ and $e$ is
called the {\it suspension disk}.  Put $\theta = 2\e/||e||$.  Then the
two angles of $\Delta$ adjacent to $e$ is at most $\arctan(2\e/||e||) <
\theta$.  Therefore $\Delta$, hence the curve $\b$, is contained in
the cones $C(\theta, \b(a), \b'(a))$ and $C(\theta, \b(b), -\b'(b))$.

Let $K$ be a polygonal knot, with edges $e_1, ..., e_n$.  A smooth
curve $\b: S^1 \to \R^3$ is an {\it allowable $\e$-approximation\/}
of $K$ if it is a union of arcs $\b_1, ..., \b_{2n}$, such that

(i)  each $\b_{2k}$ is an $\e$-allowable arc in some ball $B_k$ of
radius at most $\e$; 

(ii) each $\b_{2k-1}$ is an $\e$-suspension, such that its base arc
$E_k$ is parallel to $e_k$, and the difference between the lengths of
$e_k$ and $E_k$ is at most $\e$.

\proclaim{Lemma 2} Given any polygonal knot $K = K(v_1, ..., v_n)$ and
a regular neighborhood $N$ of $K$, there is an $\e>0$ such that any
allowable $\e$-approximation $\r$ of $K$ with the same initial point
is a knot, which is isotopic to $K$ in $N$.  \endproclaim

\proof Rescaling $\R^3$ if necessary, we may assume that the length of
each edge of $K$ is at least $3$.  Let $e_i = e(v_i, v_{i+1})$.
Denote by $m$ the minimum distance between nonadjacent edges, and by
$r$ the number given in Lemma 1.

Let $\e$ be a very small positive number (for example, $\e <
\text{min}(1, m/10n, r/10n)$).  Let $\b_1, ..., \b_{2n}$ be the arcs
of $\r$, and $B_i$ the ball containing $\b_{2i}$, as in the
definition of allowable $\e$-approximation.  Let $\Delta_i$ be the
equilateral triangles containing $\b_{2i-1}$, as in the definition of
$\e$-suspension arcs.  Denote by $v''_i, v'_{i+1}$ the initial and
ending points of $\b_{2i-1}$, respectively.  Consider the union of all
$\Delta_i$ and $B_i$.

CLAIM. {\it The triangles $\Delta_i$ are mutually disjoint, the balls
$B_i$ are mutually disjoint, and $\Delta_i$ intersects $B_j$ only if
$j=i$ or $i-1$ mod $n$, in which case they intersects at a single
point.}

Since the base arc $E_i = e(v''_i, v'_{i+1})$ of $\b_{2i-1}$ is
parallel to $e_i$ with length difference at most $\e$, and since
$d(v'_i, v''_i)$ is at most $2\e$ (the upper bound of diameters of
$B_i$), one can show by induction that $d(v_i, v'_i) \leq (3i-2)\e$,
and $d(v_i, v''_i) \leq 3i\e$.  Put $\delta = 4n\e < m/2$.  Then
$\b_{2i-1}$ is in the $\delta$-neighborhood of $e_i$, and $B_i$ is in
the $\delta$-neighborhood of $v_{i+1}$.  Since the distance between
two vertices or nonadjacent edges of $K$ is bounded below by $m$, it
follows that the balls $B_i$ are mutually disjoint, $\Delta_i$ is
disjoint from $\Delta_j$ when $i$ and $j$ are not adjacent mod $n$,
and disjoint from $B_j$ if $j$ is not equal or adjacent to $i$ mod
$n$.  Since $||E_i|| > 3 - \e$ and the height of $\Delta_i$ is at most
$\e$, the two angles of $\Delta_i$ adjacent to $E_i$ is at most
$2\e/(3-\e) < \e$.  Thus for each endpoint $v$ of $E_i$, $\Delta_i$ is
contained in a cone of angle $\e$ based at $v$ in the direction of the
tangent or negative tangent of $\b$ at $v$.  Since $\b_{2i}$ is an
$\e$-allowable arc, it follows from definition that $\Delta_i$ is
disjoint from $\Delta_{i+1}$, and they each intersects $B_i$ only at a
single point.  This completes the proof of the claim.

Since each $\b_i$ is an embedding, it follows from the claim that $\r:
S^1 \to \R^3$ is an embedding, hence is a knot.  We can isotope
$\b_{2i-1}$ via the suspension disk to the edge $E_i$.  Since
$\b_{2i}$ is unknotted in $B_i$, it can be rel $\bdd$ isotoped to a
straight arc $E'_i$ in $B_i$.  By the claim these isotopies form an
isotopy of $\r$ to the polygonal knot $K_2 = E_1 \cup E'_1 \cup
... \cup E_n \cup E'_n = K(v'_1, v''_1, ..., v'_n, v''_n)$.  Since
$d(v'_i, v''_i)$ is very small, by Lemma 1(ii) $K_2$ is isotopic to the
knot $K(v''_1, ..., v''_n)$, which is isotopic to $K$ by Lemma 1(iii).
\endproof

Let $A$ be a compact 1-manifold.  A smooth map $\a: A \to S^2$ is {\it
admissible\/} if (i) $\a$ is an embedding, and (ii) it has no
antipodal points, i.e., $\a(t) \neq - \a(s)$ for all $t \neq s$.
Denote by $\eta: S^2 \to S^2$ the antipodal map, and by $\rho: S^2 \to
P^2$ the standard double covering map onto the projective plane $P^2$.
Then $\a$ is admissible if and only if $\rho \circ \a : A \to P^2$ is
a smooth embedding.

\proclaim{Lemma 3} Suppose $Y$ is the disjoint union of finitely many
circles, and suppose $A$ is a compact submanifold of $Y$.  Let $\a: A
\to S^2$ be an admissible map such that each circle component of
$\a(A)$ bounds a disk $\Delta$ with interior disjoint from $\a(A)$ and
$\eta(\Delta)$.  Then $\a$ extends to an admissible map $\hat \a: Y
\to S^2$.  \endproclaim

\proof Let $I$ be the closures of components of $Y-A$.  We need to
extend $\a$ to an admissible map $\hat \a : A \cup I \to S^2$ which
still satisfies the assumption of the lemma.  The result would then
follow by induction.  If $I$ is a circle, define $\hat \a : I \to S^2$
to be a smooth map embedding $I$ into a small disk $D$ of $S^2$ such
that $D$, $\rho(D)$ and $\a(A)$ are mutually disjoint.  So suppose $I$
is an interval with endpoints $u_1, u_2$ on a component $Y_0$ of $Y$.
Denote by $\tilde \a = \rho \circ \a$.

If $J = Y_0-\Int I$ is connected, then by assumption $\tilde \a$ is an
embedding, so there is a small disk neighborhood $D$ of $\tilde \a(J)$
which is disjoint from $\tilde \a(A-J)$.  Let $D_1$ be the component
of $\rho^{-1}(D)$ containing $\a(J)$, and extend $\a$ to a smooth
embedding $\hat \a : A \cup I \to S^2$ so that $\hat \a (I) \subset D$.

Now suppose $J$ is disconnected.  Let $J_1, J_2$ be the components of
$J$ containing $u_1, u_2$ respectively.  Let $K_1, ..., K_r$ be the
circle components of $A$, and let $D_i$ be the disk on $P^2$ bounded
by $K_i$.  By assumption $\tilde \a(J_i)$ are in $P^2 - \cup D_i$, so
there are two non-homotopic arcs $\tilde \r_1, \tilde \r_2 : I \to
P^2$ such that $\tilde \r_i \cup \tilde \a : I \cup A \to P^2$ is a
smooth embedding.  One of the $\tilde \r_i$ lifts to a path $\r: I \to
S^2$ connecting $u_1$ to $u_2$.  It follows that $\r \cup \a : I \cup
A \to S^2$ is the required extension.  \endproof

\proclaim{Lemma 4} Suppose $\a = (\a_1, \a_2, \a_3): I= [a,b] \to S^2$
is an admissible curve intersecting $S_1$ transversely at two points
in the interior, and $\a_3(a) > 0$.  Then there is a function $f: I
\to \R_+$ such that (i) $f(t) = 1$ in a neighborhood of $\bdd I$, and
(ii) the integral curve $\b$ is unknotted in $Z[z_1,z_2]$, where $z_1 =
\b_3(a)$ and $z_2 = \b_3(b)$.  \endproclaim

\proof By assumption $\a_3$ has exactly two zeroes $u, v\in I$,
($u<v$), so $a_3(t) < 0$ if and only if $t \in (u,v)$.  Since $\a$ is
admissible, $\a(u) \neq \pm \a(v)$, so by a rotation along the
$z$-axis if necessary we may assume that $\a_1(u),\a_1(v) > 0$, and
$\a_2(u), \a_2(v)$ have different signs.  Without loss of generality
we may assume that $\a_1(t), \a_2(t) > 0$ when $t$ is in an
$\e$-neighborhood of $u$, and $\a_1(t) > 0$, $\a_2(t) < 0$ when $t$ is
in an $\e$-neighborhood of $v$, where $0<\e< \text{min} (u-a, b-v)$.

We start with the constant function $f(t) = 1$ on $I$, and proceed to
modify $f(t)$ so that $f(t)$ and the integral curve $\b = \b(f, \a, t_0,
v_0)$ satisfy the conclusion of the lemma.  Put $\b = (\b_1(t),
\b_2(t), \b_3(t))$, and choose the base point $v_0$ so that $\b(u) = 0$.
Thus
$$ \b_i(t) = \int_u ^t f(t) \a_i(t) \, dt.$$ Since $\a_1(u), \a_2(u) >
0$, and $\b_1(u) = \b_2(u) = 0$, by enlarging $f(t)$ in a small
$\e$-neighborhood of $u$, we may assume that $\b_1(t), \b_2(t) > 0$
for all $t \in (u,v]$.  Since $\a_2(t) < 0$ in a neighborhood of $v$,
we may then enlarge $f(t)$ near $v$ so that $\b_2(v) = \b_2(u) = 0$.
This does not affect the fact that $\b_1(t) > 0$ for $t \in (u, v]$,
and $\b_2(t) > 0$ for $t \in (u, v)$.

The function $\b_3$ is descending in $[u,v]$ because $\a_3(t)$ is
negative in this interval.  Thus $\b_3(v) < \b_3(u)$.  Since $\a_3$ is
positive in $[a,u]$ and $[v,b]$, $\b_3$ is increasing in these
intervals.  We may now enlarge $f(t)$ in $(u-\e, u)$ and $(v, v+\e)$,
so that $z_3 = \b_3(u-\e) < \b_3(v)$ and $z_4 = \b_3(v+\e) > b_3(u)$.
Thus the curve $\b$ on $[u-\e, v+\e]$ is a proper arc in $Z[z_3,
z_4]$.  We want to show that it is unknotted.

By the above, the curve $\b[u,v]$ lies in $Z[z_3,z_4] \cap
Y[0,\infty)$, with endpoints on the $xz$-plane.  Since $\b_3$ is
descending on $[u,v]$, $\b$ is rel $\bdd$ isotopic in $Z[z_3,z_4] \cap
Y[0,\infty)$ to a straight arc $\hat \b[u,v]$ on the $xz$-plane.
Since $\a_2(t) > 0 $ for $t \in [u-\e, u]$, and $\b_2(u) =0$, we have
$\b_2(t) < 0$ for $t \in [u-\e, u]$.  Similarly, since $\a_2(t) < 0$
near $v$, we have $\b_2(t) < 0$ for $t \in [v, v+\e]$.  Therefore, the
above isotopy is disjoint from the arcs $\b[u-\e, u]$ and $\b[v,
v+\e]$, hence extends trivially to an isotopy of $\b[u-\e, v+\e]$,
deforming $\b[u-\e, v+\e]$ to the curve $\hat \b = \b[u-\e, u] \cup
\hat \b[u,v] \cup \b[v, v+\e]$.

Since $\a_1(t)$ is positive near $u, v$, $\b_1$ is increasing in
$[u-\e, u]$ and $[v, v+\e]$.  Since $\hat \b$ is a straight arc
connecting $\b(u)$ and $\b(v)$, and $\b_1(v) > \b_1(u)$ by the above,
the first coordinate function of $\hat \b$ is also increasing in $[u,
v]$.  It follows that the first coordinate of $\hat \b$ is increasing
in $[u-\e, v+\e]$, therefore, $\hat \b$ is unknotted in $Z[z_3,z_4]$,
hence is rel $\bdd$ isotopic to a straight arc $\tilde \b$ in
$Z[z_3,z_4]$.

Since $\b_3(t)$ is increasing on $[a, u-\e] \cup [v+\e, b]$, the above
isotopy extends trivially to an isotopy deforming $\b: I \to \R^3$ to
the curve $\b[a,u-\e] \cup \tilde \b \cup \b[v+\e, b]$.  Since the
third coordinate of this curve is always increasing, it is unknotted
in $Z[z_1,z_2]$, where $z_1 = \b_3(a)$ and $z_2 = \b_3(b)$.
Therefore, $\b$ is also unknotted in $Z[z_1,z_2]$.  \endproof

Given $a \in \R$ and $\d > 0$, let $\varphi = \varphi[a,\d](x)$ be a
smooth function on $\R^1$ which is symmetric about $a$, $\varphi(a) =
1$, $\varphi(x) = 0$ for $|x - a| \geq \d$, and $0\leq \varphi(x) \leq
1$ for all $x$.  Given $a,b\in \R$ with $a < b$, let $\psi(x) =
\psi[a,b](x)$ be a smooth monotonic function such that $\psi(x) = 0$
for $x\leq a$, and $\psi(x) = 1$ for $x \geq b$.  Such functions
exist, see for example [GP, Page 7].

For any point $p \in S^2$, denote by $U(p, \e)$ the
$\e$-neighborhood of $p$ on $S^2$, measured in spherical
distance.  Thus for any $q\in U(p, \e)$, the angle between $p, q$
(considered as vectors in $\R^3$) is less than $\e$.  

\proclaim{Lemma 5} Let $0< \e < \pi/8$, and let $\a = (\a_1,\a_2,\a_3)
: I = [a_1, a_2] \to S^2$ be an admissible arc transverse to $S_1$,
such that $\a_3(a_i)>\e$.  Let $\mu > 0$.  Then there is a smooth
positive function $f(t)$ such that (i) $f(t) = 1$ near $a_i$, and (ii)
the intergral curve $\b = \b(f, \a, t_0, v_0)$ is an $\e$-allowable
arc in a ball of radius $\mu$ in $\R^3$.  \endproclaim

\proof  Notice
that $U(\a(a_i), \e)$ are on the upper half sphere $S^2_+$.
Choose $0.1 > \d > 0$ sufficiently small, so that $\a(t) \in
U(\a(a_i), \e)$ for $t$ in a $\d$-neighborhood of $a_i$.  Choose
$c_0 = a_1 + \d, c_1, ..., c_p = a_2 - \d$ so that the curve $\a(I_j)$
intersects $S_1$ exactly twice in the interior of $I_j = [c_{j-1},
c_j]$, $j=1,...,p$.

By Lemma 4  applied to each $I_j$, we see that there is a function
$f_1(t)$ on $I$, such that $f_1(t) = 1$ near $c_i$ and on
$[a_1,c_0]\cup[c_p,a_2]$, and the part $\b_1[c_0, c_p]$ of the
integral curve $\b_1 = \b(f_1, \a, t_0, v_0)$ is unknotted in $Z[z_0,
z_p]$, where $z_i = \b(c_i)$.  Without loss of generality we may
choose $t_0 = c_0$ and $v_0 = 0$.  Since the curve is compact, the
isotopy is within a ball, so there is a disk $D$ in $\R^2$, such that
$\b_1[c_0, c_n]$ is unknotted in $D \times [z_0,z_p]$.  Choose $N$
large enough, so that the ball $B(N)$ of radius $N$ centered at the
origin contains both $D\times [z_0,z_p]$ and the curve $\b_1$ in its
interior.  We want to modify $f_1(t)$ on $[a_1,c_0)\cup (c_p,a_2]$ to
a function $f_3(t)$, so that $\b_3=\b(f_3,\a,c_0,0)$ is an
$\e$-admissible curve in $B(10N)$, and $f_3(t) = 10N/\mu$ near
$\bdd I$.

First, consider the function $$f_2 = f_1 + (\frac{10N}{\mu}-1) (1-
\psi[a_1, a_1+\e_1] + \psi[a_2-\e_1, a_2]),$$ where $\e_1$ is a very
small positive number, say $\e_1 < \text {min}(\d, \mu/10)$.  By the
property of the $\psi$ functions, we have $f_2(t) = f_1(t)$ for $t \in
[c_0, c_p]$, and $f_2(t) = 10N/\mu$ near $a, b$.  Let $\b_2 = \b(f_2,
\a, z_0, 0)$.  Since $\e_1$ is very small, one can show that
$||\b_2(t)|| < 2N$ for all $t \in [a,b]$.
Let $b_1, b_2$ be positive real numbers.  Define
$$ f_3(t) = f_2(t) + b_1\, \varphi[a_1+\delta/2, \delta/4](t) +
b_2 \,\varphi[a_2-\delta/2, \delta/4](t).$$
Let $ \b_3$ be the integral curve $\b_3(f_3,\a,z_0, 0)$.  
Since $L$ is a polygonal knot, $\sum e_j = 0$, so we have
$$\b_3(a_2) = \b_2(a_2) + b_2 \int_{z_p}^b \varphi[a_2-\delta/2,
\delta/4](t) \a(t) \, dt = \b_2(a_2) + b_2 v_2.$$ Since $\a(t) \in
U(\a(a_2), \e)$ and $\e < \pi/8$, the vector $v_2$ above is nonzero.
Since $||\b_2(a_2)|| < 2N$, we may choose $b_2>0$ so that
$||\b_3(a_2)|| = 10N$.  Similarly, choose $b_1>0$ so that $\b_3(a_1) =
10N$.

Consider a point $t \in [c_p, a_2]$ such that $||\b_3(t)|| \geq 10N$.
Let $\theta(t)$ be the angle between $\b_3(t)$ and $\b_3'(t)$.  Put
$u_0 = \b_3(c_p)$, and notice that $||u_0|| < N$.  Since $\a(t) \in
U(\a(a_2), \e)$, the curve $\b_3[c_p, b]$ lies in the cone $C(\e, u_0,
\a(a_2))$, so the angle between $(\b_3(t) - u_0)$ and $\a(t)$ is at
most $2\e$.  We have
$$ \align
\cos \theta(t) & = \frac {\b_3(t) \cdot \a(t)}{||\b_3(t)||} 
= \frac {(\b_3(t) - u_0) \cdot \a(t) + u_0 \cdot \a(t)}{||\b_3(t)||}
\\
& \geq \frac {(10N - N) \cos (2\e) - N}{10N}  =  0.9 \cos(2\e)
- 0.1 \\
& > \frac 12
\endalign
$$
Therefore, $\theta(t) < \pi/3$.  In particular, this implies that the
norm of $\b_3(t)$ is increasing if it is at least $10N$ and $t \in
[c_p, a_2]$; but since $||\b_3(a_2)|| = 10N$, it follows that $\b_3(t)
\in B(10N)$ for $t \in [c_p, a_2]$.  Similarly, one can show that this
is true for $t\in [a_1, c_0]$.  Therefore, $\b_3$ is a proper arc in
$B(10N)$.  It is unknotted because its third coordinate is increasing
on $[a_1, c_0] \cup [c_p, a_2]$ and the curve $\b_3[c_0, c_p] =
\b_1[c_0, c_p]$ is unknotted in $D \times [z_0, z_p]$, with
$\b_3(c_0)$ on $D\times z_0$.  

We need to show that the cone $C(\e, \b_3(a_2), \b_3'(a_2))$
intersects $B(10N)$ only at the cone point, but this is true because
$\e + \theta(a_2) < \pi/8 + \pi/3 < \pi/2$.  Similarly for $C(\e,
\b_3(a_1), -\b'_3(a_1))$.  Also, notice that the cone $C(\e,
\b_3(a_2), \b'_3(a_2))$ lies above the $xy$-plane, while $C(\e,
\b_3(a_1), -\b'_3(a_1))$ lies below the $xy$-plane, so they are
disjoint.  It follows that $\b_3$ is an $\e$-allowable curve in
$B(10N)$.

Finally, rescale the curve by defining $f(t) = f_3(t) \mu/10N$, and
$\b = \b(f, \a, c_0, 0)$.  Then $\b$ is an $\e$-allowable curve in a
ball of radius $\mu$, and $f(t) = 1$ near $\bdd I$.  \endproof

\proclaim{Lemma 6} Suppose the integral curve $\b = \b(f, \a, a, 0)$
is a round $\e$-suspension.  Then for any $k \in [\frac 12, \frac
32]$, there is a positive function $g(t)$ such that (i) $g(t) = f(t)$
near $a,b$, and (ii) the integral curve $\r = \b(g,\a,a,0)$ is a
$(k\e)$-suspension with $\r(b) -\r(a) = k (\b(b) - \b(a))$.
\endproclaim

\proof Without loss of generality we may assume $[a,b] = [-1,1]$.  Set
up the coordinate system so that $\b$ lies in the triangle with
vertices $\b(a) = (0,0,0)$, $\b(b) = (2u,0,0)$ and $(u,\e,0)$, where
$2u = ||\b(b) - \b(a)||$.  Put $\a = (\a_1, \a_2, \a_3)$.  Then
$\a_3(t) = 0$, and $\a_2(-t) = -\a_2(t)$.  Consider the smooth
function $\phi = \psi[-1+\d, -1+2\d] - \psi[1-2\d, 1-\d]$.  It is an
even function with $\phi(t) = 1$ when $|t| \leq 1-2\d$, and $\phi(t) =
0$ when $|t| \geq 1-\d$.  Let $g(t) = c + p\phi(t)$, where $p > -c$ is
a constant.  Since $|\phi(t)| \leq 1$, $g(t)$ is a positive function.
We have
$$\r(1) = \int_{-1}^1 g(t) \a(t) \, dt = \b(1) + p \int_{-1}^1 \phi(t)
\a(t) \, dt.$$ Since $\phi(t)$ is even and $\a_2$ is odd, $\r_2(1) =
\r_3(1) = 0$.  When $\d$ approaches $0$, the integral
$$ c \int_{-1}^1 \phi(t) \a_1(t) \, dt$$ approaches $\b_1(1) = 2u$.
Hence for any $s \in [u, 3u]$, we may choose $\d$ small and $p \in
(-c, c)$ so that $\r(1) = (s, 0, 0)$.  Note that $\r'_1(t) = (c - p)
\a_1(t) > 0$, so $\r$ is an embedding.

Consider $\r$ and $\b$ as curves on the $xy$-plane.  Then The tangent
of $\r$ at $t$ is given by 
$$\frac {\r'_2(t)}{\r'_1(t)} = \frac {g(t)\a_2(t)}{g(t)\a_1(t)} =
\frac {\a_2(t)}{\a_1(t)},$$ which is the same as the tangent slope of
$\b$ at $t$, and hence is bounded above by $\e$.  Thus $\r$ is below
the line $y = \e x$ on the $xy$-plane.  Similarly, it is below the
line $y = -\e(x-\r_1(1))$.  It follows that $\r$ is a
$k\e$-suspension, where $k = \r_1(1) / \b_1(1) < 2$.  \endproof

\proclaim{Theorem 7} Given any tame link $L$ in $S^3$ and any
neighborhood $\eta(L)$ of $L$, there is a smooth link in $\eta(L)$ which is
isotopic to $L$ in $\eta(L)$, and has no parallel or antiparallel
tangents.  \endproclaim

\proof Without loss of generality we may assume that $L = K_1 \cup
... \cup K_r$ is an oriented polygonal link in general position, with
oriented edges $e_1, ..., e_m$, which are also considered as vectors
in $\R^3$.  Let $d$ be the minimum distance between nonadjacent edges.
We may assume each $K_i$ has at least four edges, so $d$ is also an
upper bound on the length of $e_i$.  For any $\e_1 \in (0, d/3)$ the
$\e_1$-neighborhoods of $K_i$, denoted by $N(K_i)$, are mutually
disjoint.  Choosing $\e_1$ small enough, we may assume that $N(K_i)$
are contained in $\eta(L)$.  By Lemma 2 there is an $\e>0$, such that
any allowable $\e$-approximation of $K_i$ is contained in $N(K_i)$ and
is isotopic to $K_i$ in $N(K_i)$.  Note that $\e \leq d/3$.  We will
construct such an approximation $\hat K_i$ for each $K_i$, with the
property that $\hat L = \hat K_1 \cup ... \cup \hat K_r$ has no
parallel or antiparallel tangents.  Since $N(K_i)$ are mutually
disjoint, the union of the isotopies from $\hat K_i$ to $K_i$ will be
an isotopy from $\hat L$ to $L$ in $N(L)$.

Consider the unit tangent vector of $e_i$ as a point $p_i$ on $S^2$,
which projects to $\hat p_i$ on $P^2$.  Since $L$ is in general
position, $\hat p_1, ..., \hat p_n$ are mutually distinct, so by
choosing $\e$ smaller if necessary we may assume that they have
mutually disjoint $\e$-neighborhoods $\hat D_1, ..., \hat D_n$, which
then lifts to $\e$-neighborhoods $D_1, ..., D_n$ of $p_1, ..., p_n$.
Adding some edges near vertices of $L$ if necessary, we may assume
that the angle between the unit tangent vectors of two adjacent edges
$e_i, e_{i+1}$ of $L$ (i.e.\ the spherical distance between $p_i$ and
$p_{i+1}$), is small (say $< \pi/2$).

Bend each edge $e_j$ a little bit to obtain a round
$(\e/2)$-suspension $\hat e_j : I_j \to \R^3$ with $||\hat e'_j(t)|| =
1$ (so the length of $I_j$ equals the length of the curve $\hat e_j$).
Then its derivative $\hat e'_j$ is a map $I_j \to S^2$ with image in
$D_j$ because $D_j$ has radius $\e$.  Let $Y = \cup S^1_i$ be a
disjoint union of $r$ copies of $S^1$, and let $A = \cup I_j$ be the
disjoint union of $I_j$.  Embed
$A$ into $Y$ by a map $\eta$ according to the order of $e_i$ in $L$.
More precisely, if $e_j$ and $e_k$ are edges of $L$ such that the
ending point of $e_j$ equals the initial point of $e_k$ then the
ending point of $\eta(I_j)$ and the initial point of $\eta(I_k)$
cobounds a component of $Y-\eta(A)$.  The union of the maps $\hat e'_j
\circ \eta^{-1}$ defines a map $\eta(A) \to S^2$, which is admissible
because the disks $\hat D_j$ on $P^2$ are mutually disjoint.  By Lemma
3, it extends to an admissible map $\hat \a: Y \to S^2$.  It now
suffices to show that each $K_i$ has an allowable $\e$-approximation
$\hat K_i: S^1_i \to \R^3$, with $\hat \a |_{ S_i^1}$ as its unit
tangent map.

The construction of $\hat K_i$ is independent of the other components,
so for simplicity we may assume that $L = K(v_1, ..., v_n)$ is a knot,
with edges $e_i = e_i(v_i, v_{i+1})$.  Since $L$ is in general
position, the three unit vectors $p_1, p_2, p_3$ of the edges $e_1,
e_2, e_3$ are linearly independent, so there is a positive number
$\delta < \e$, such that the ball of radius $\delta$ centered at the
origin is contained in the set $\{ \sum u_i p_i \, | \, \e > |u_i|
\}.$

We may assume that the intevals $I_j = [a_j, b_j]$ are sub-intervals
of $I = [0, a_{n+1}]$, with $a_1 =0$ and $b_j < a_{j+1}$.
Put $\hat I_j = [b_j, a_{j+1}]$.  Without loss of generality we may assume
that the function $\eta: A \to S^1$ defined above is the restriction
of the function $\eta : I \to S^1$ defined by $\eta(t) = \text{exp}(2\pi i
t/a_{2n})$.  Put $\a = \hat \a \circ \eta : I \to S^2$.
Thus $\a(t) = \hat e'_j(t)$ on $I_j$.

Consider the restriction of $\a$ on $\hat I_j = [b_j, a_{j+1}]$.  We have
assumed that the spherical distance between $p_j$ and $p_{j+1}$ is at most $\pi/2$.  Since 
$\a(b_j) \in D_j$ and $\a(b_{j+1}) \in D_{j+1}$, 
the spherical distance between
$\a(b_j)$ and $\a(a_{j+1})$ is at most $\pi/2 + 2\e$.  As $\e$ is
very small, we may choose a coordinate system for $\R^3$ so that the
third coordinate $\a_3$ of $\a$ is greater than $\e$ at $b_j$ and
$a_{j+1}$, and by transversality theorem we may further assume
that $\a$ is transverse to the circle $S_1 = S^2 \cap \R_{xy}$ in this
coordinate system.  Now we can apply Lemma 5 to get a function
$f_j(t)$ on $\hat I_j$ such that $f_j(t) = 1$ near $\bdd \hat I_j$,
and the integral curve $\r_j = \b(f_j, \a|_{\hat I_j})$ is an
$\e$-allowable curve in a ball of radius $\d/2n$.  Extend these
$f_j$ to a smooth map on $I$ by defining $f(t) = 1$ on $I_j$.

Consider the integral curve $\b = \b(f, \a)$.  It is the union of $2n$
curves $\b_i, \hat \b_i$ defined on $I_i$ and $\hat I_i$, where $\b_i =
\b(f|_{I_i}, \a|_{I_i}, a_i, \b(a_i))$ is a translation of $\hat e_j$ 
because $f|_{I_i} = 1$;
and $\hat \b_i = \b(f|_{\hat I_i}, \a|_{\hat I_i}, b_i, \b(b_i))$ is a
translation of $\r_i$ because $f|_{\hat I_i} = f_i$.  We have
$$\align ||\b(a_{n+1}) - \b(a_{0})|| & \leq \sum_1^n ||\b(a_{j+1}) -
\b(b_j)|| + ||\sum _1^n (\b(b_j) - \b(a_j))|| \\ & \leq \sum 2(\d/2n)
+ ||\sum e_j|| = \d. \endalign$$ By the definition of $\d$, there are
numbers $u_i \in [-\e, \e]$, such that $\b(2n) - \b(0) = \sum_{i=1}^3
u_i p_i$.  Notice that $|u_i| < \e < ||e_j|| / 2$, so by Lemma 6, we
can modify $f(t)$ on $[a_j+\e_1, \b_j-\e_1]$ for $j=1,2,3$ and some
$\e_1>0$, to a function $g(t)$, so that the integral curve $\r = \b(g,
\a)$ on $I_j$ is an $\e$-suspension with base arc the vector $e_j +
u_j p_j$.  Now we have $\r(a_{n+1}) = \r(0)$, so $\r$ is a closed
curve.  Since $\r'(t) = \a(t)$ near $0$ and $2n$ and $\a$ induces a
smooth map $\hat \a : S^1 \to S^2$, it follows that $\r$ induces a
smooth map $\hat \r: S^1 \to \R^3$.  

From the definition we see that $\hat \r$ is an allowable
$\e$-approximation of $K$.  This completes the proof of the theorem.
\endproof

ACKNOWLEDGEMENT.  I would like to thank Colin Adams for
raising the problem, and for some helpful conversations.

\enddocument